\documentclass{amsart}

\newtheorem*{theorem}{\rm THEOREM}
\newtheorem{lem}{\rm LEMMA}

 \theoremstyle{remark}
\newtheorem*{rem}{\rm REMARK}

 \numberwithin{equation}{section}

\def\R{{\mathbb R}}   % reals
\def\Q{{\mathbb Q}}   % rationals
\def\di{{\triangle}}   %discriminant
\def\e{{\mathrm{e}}}

\newcommand{\pfrac}[2]{\genfrac{(}{)}{}{}{#1}{#2}} %proper frac
\newcommand{\E}{\zeta_{Q} } %Epstein's zeta
\newcommand{\half}{\frac{1}{2}}
\newcommand{\cc}{8}
\newcommand{\dd}{1/2}
\begin{document}

\title{On the zeros of the Epstein zeta function}
\author{Anirban Mukhopadhyay, Krishnan Rajkumar, Kotyada Srinivas}
\address{Institute of Mathematical Sciences,
CIT Campus, Tharamani, Chennai 600 113, India}
\email[Anirban Mukhopadhyay]{anirban@imsc.res.in}
\email[Kotyada Srinivas]{srini@imsc.res.in}
\email[Krishnan Rajkumar]{rkrishnan@imsc.res.in}

\begin{abstract} 
In this article, we count the number of consecutive zeros of the Epstein zeta-function,
associated to a certain quadratic form, on the critical line with ordinates lying in $ [0, T], T $ sufficiently large and which are separated apart by a given positive number $V$.
\end{abstract}

 \maketitle

\section{Introduction}{\label{intro}
The study of the distribution of zeros of the zeta-functions on the {\it critical line} is a 
fundamental problematique in analytic number theory. A related problem is to find a good upper
bound for the difference beteween consecutive zeros on the critical line. Such a study natually
leads one to ask the following 

\medskip
\noindent {\bf Question 1.} \textit{Given 
$T$ sufficiently large, find $H = H(T)$ such that the interval $ [ T , T+H ] $ contains the ordinate of a 
 zero on the critical line of the corresponding zeta function under consideration? }

\medskip

\noindent {\bf Question 2.} \textit{Given $T$ sufficiently large and $V>0$, how many  consecutive zeros of 
the zeta-function under consideration are there on the critical line with ordinates in $[0,T]$ which are atleast $V$ distance apart?}

\medskip
Question 1 has been studied in great detail from the time of Hardy and Littlewood. 
We shall briefly mention
few known results below. The main object of this note is to  study the second 
question with respect to the Epstein zeta function
$ \zeta_Q(s) $ (see \eqref{int0e} ). We shall address this issue in section \ref{result}.

In a classical paper \cite{H-L},   Hardy and Littlewood
gave an answer to Question 1 for the  Riemann zeta function, showing that one can take 
$H = T^{1/4 + \varepsilon} $. 
This was subsequently improved to $H = T^{1/6} (\log T)^{5+\varepsilon}  $ 
by Moser \cite{moser}  and 
to $ H = T^{1/6 + \varepsilon}  $ by  Balasubramanian  \cite{B} ,
and the latest result is by Karatsuba \cite{kara} who showed that $H= T^{5/32 } (\log T)^2 $ holds. 
It must be mentioned here that Ivi\'c ( page 261, \cite{ivic} ) improved Karatsuba's exponent 
$5/32 = 0.15625$ to $ 0.1559458\dots  $.  
Following the method of Hardy and Littlewood, Potter and Titchmarsh obtained analogous 
result for the Epstein zeta-function, which we shall discuss in section \ref{result}. 
For the zeta function, $\zeta_K(s,C) $ of an ideal class $C$ in a 
quadratic number field $K$, Bruce Bernt \cite{bruce} 
showed that $ H = T^{1/2 +\varepsilon }  $ holds. This was improved to 
$ H= T^{1/2} \log T  $ by Sankaranarayanan \cite{San}.
In \cite{jutila1}, Jutila developed a method of transforming certain exponential 
sums into another sum which is much easier to handle ( see section \ref{summation} ) 
and thereby  proved that  $H= T^{1/3 + \varepsilon}$ holds for the $L$-function 
associated to a cusp form \eqref{cusp} for the full modular group.
 
\medskip

A well-known theorem of Selberg \cite{Sel} states that 
the function $\zeta(\half+it)$
has at least $\gg T\log T$ zeros in the interval $[0,T]$, and this is best possible
upto a value of the implied constant. In fact, Heath-Brown \cite{heath-brown} 
showed that the same estimate holds 
even if we count only simple zeros. Hence the average gap between critical zeros of $ \zeta(s)$ is 
$(\log T)^{-1}$. It is in this context that the Question 2 becomes relevant!

Let $R_1:=R_1(V)$ be the number of gaps  between consecutive zeros of 
$ \zeta(\half+it)$, 
with ordinate between $0$ and $T$,  which are larger than $V$. Then trivially $R_1\ll TV^{-1}$.
Karatsuba \cite{karatsuba} showed that 
$R_1\ll TV^{-3/2}$ for $V=T^{\delta}$ for any $\delta>0$.
Ivi\'{c} and Jutila \cite{ivic-jutila} made substantial improvement in the
case of $ \zeta(s)$ as well as the $L$-function associated to a certain cusp form. 
We state their results below as Theorem A and B:\\

\noindent THEOREM A:~~
Let $R_1$ be as defined above. Then uniformly
\begin{equation}\label{thmA1}
R_1 \ll T V^{-2}\log T
\end{equation}
and 
\begin{equation}\label{thmA2}
R_1 \ll T V^{-3}\log^5 T.
\end{equation}

For $V\ll \log T$ these results are trivial and also $R_1=0$ for $V=T^{1/6}$
which follows from known results mentioned earlier on the gap between 
critical zeros for $ \zeta(s)$.
Thus, the result is non-trivial in the range
$\log T\ll V \le T^{1/6-\epsilon}.$ In this article we shall be always
considering $V$ in the range $T^{\epsilon}\ll V \le T^{1/6-\epsilon}$.

\medskip

Let $a(n)$'s denote the Fourier coefficients
of a cusp form $f$ of weight $\kappa$ for the full modular group. The
$L-$series associated with $f$ is given by the  Dirichlet series
\begin{equation}\label{cusp}
\varphi (s)=\sum_{n=1}^{\infty} \frac{a(n)}{n^s}, \qquad Re(s)>(\kappa+1)/2 .
\end{equation}

Let $R_2:=R_2(V)$ denote the number of gaps of length at least $V$ between consecutive
zeros of $\varphi(\kappa/2 + it)$ in the interval $[0,T]$.

\medskip

\noindent{THEOREM B:} Suppose that the $a(n)$'s defined above are real and let 
$ V\gg \log^5 T $.
Then
\begin{equation}\label{thmB1}
R_2 \ll T V^{-2} \log T
\end{equation}
and 
\begin{equation}\label{thmB2}
R_2 \ll T^2 V^{-6} \log^6 T.
\end{equation}

\medskip

In this note, we obtain an analogue of the results \eqref{thmA1} and \eqref{thmB1} for the 
Epstein zeta-function, $\zeta_Q(s)$. We state the main result in the
next section after a brief survey of the known results on the zeros of $\zeta_Q(s)$.
In section \ref{lemmas} and \ref{summation} we shall state some Lemmas, section \ref{proof}
will deal with the proof of the main result.  

\section{ Main result}\label{result}

Let $Q(x,y) = a x^2+b x y + c y^2$ be a binary positive-definite integral quadratic form. 
Let $r_Q(n)$ be the number of solutions of $Q(x,y)=n$. The Epstein zeta-function associated to 
$Q(x,y)$ is defined as 
\begin{equation}\label{int0e}
\zeta_Q(s)=\sum_{(x,y)\neq(0,0)} \frac{1}{Q(x,y)^s} = \sum_{n=1}^{\infty} \frac{r_Q(n)}{n^s},
 \qquad \sigma > 1. 
\end{equation}

It is well known that $\zeta_Q(s)$ can be analytically continued to the entire complex plane except
for $s=1$ where it has a simple pole with residue $2\pi/{\sqrt{\Delta}}$
where $\triangle=4ac-b^2 > 0$.
It also has the 
functional equation
\begin{equation}\label{int1e}
\left(\frac{\sqrt{\triangle}}{2\pi}\right) ^s \Gamma(s)\E(s)
={\left(\frac{\sqrt{\triangle}}{2\pi}\right)}^{1-s}\Gamma(1-s)\E(1-s).
\end{equation}
We suppose
throughout the paper that $\Delta $ is not a square, so that
$\sqrt{\triangle}$ is irrational; 
other cases like $\triangle = 4$ related to the form
 $Q(x,y)=x^2+y^2$ are either easier or well-known.

Let $h(-\triangle)$ be the number of inequivalent forms with discriminant $-\triangle$. 
It is well known that there is a one-one correspondence between the equivalence classes 
of quadratic forms of discriminant $-\triangle$ 
and the ideal classes of $K=\Q(\sqrt{-\triangle})$.
In fact, the Dedekind zeta-function $\zeta_K(s)$ for the number field $K$ can be written as
$$
\zeta_K(s) = \frac{1}{l} \sum_r \zeta_{Q_r}(s),
$$
where $Q_r$ runs over the distinct equivalence classes and $l$ is the number of units in K.

The Dedekind zeta-function has an Euler
product and it belongs to that class of Dirichlet series for which
Riemann's hypothesis can be reasonably expected. In contrast, if
the class number $h(-\triangle)$ of the quadratic forms with
discriminant $-\triangle$ exceeds one, then $\E(s)$ has no Euler
product and the Riemann hypothesis fails to hold for it. 
Indeed, when $a,b,c$ are integers, $-\triangle$ is a fundamental discriminant
and the class number $h(-\triangle)>1$, Davenport and Heilbronn \cite{daven-heil}
had shown that $\E(s)$ has infinitely many zeros in the half plane $\sigma>1$ 
aribitrarily close to the line $\sigma=1$.
S. M. Voronin \cite{Vor} proved that if $h(-\triangle)>1$, then the
number of zeros of $\E(s)$ in the rectangle
$\sigma_1\leq\sigma\leq\sigma_2$, $|t|\leq T$, with $1/2< \sigma_1
<\sigma_2 \leq 1$ and $T$ sufficiently large, exceeds $c(\sigma_1
,\sigma_2 ) T$, where $c(\sigma_1 ,\sigma_2 )>0$.
Chowla and Selberg \cite{chow-sel} showed the existence of a real zero $s$ with
$\half < s < 1$ when $a=1,b=0$ and $c$ is large enough. 
Bateman and Grosswald \cite{bate-gross} improved it by showing that $\E(s)$ has a real zero
between $\half$ and $1$ if $k=(\sqrt{\triangle})/2a>7.0556.$.

Stark \cite{stark} gave a striking description about 
the distribution of zeros of $\E(s)$ in a 
bounded region. More precisely, he proved that there 
exists a number $K$ such that if $k>K$, 
then all the zeros of $\E(s)$ in the region
$-1 < \sigma < 2$, $-2k \le  t \le 2k$ are simple zeros and with the exception of two real
zeros between $0$ and $1$, all are on the line $\sigma=\half$.
Moreover, in the same paper he proved the following zero-density result:
Let $N(T,Q)$ denote the number of zeros of $\E(s)$ in the region
$-1 < \sigma < 2$, $0 \le  t \le T$.
If $k>K$ and $0 < T \le 2k$, then 
$$ N(T,Q)=(T/{\pi})\log (kT/{\pi e})+O(h(T+3)),$$
where $h(x)=(\log x)^{1/3} (\log\log x)^{1/6}$ and the constant
in $"O"$ does not depend on $k$.

Though $\E(s)$, in general, has thus infinitely many zeros \emph{off} the
critical line $\sigma = 1/2$,
it has nevertheless infinitely many zeros \emph{on}
the critical line, and in fact it was proved by Potter
and Titchmarsh \cite{P-T} that there is a zero $1/2 +
i\gamma$ of $\E$ with $\gamma \in [\,T,\, T+ T^{1/2 +
\varepsilon}\,]$, for any fixed $\varepsilon >0$ and $T\geq
T(Q, \varepsilon)$. Sankaranarayanan \cite{San} sharpened this
by showing the same for intervals of the type 
$[\, T,\, T+cT^{1/2} \log T\,]$. 
For a long time improving the exponent of $T$
below $1/2$ was considered to be a challenging problem. First result in this
direction is by Jutila and Srinivas \cite{srini-jut}, where they prove 
that for any positive definite binary integral
quadratic form $Q$ and for any fixed $\varepsilon >0$ and 
$T\geq T(\varepsilon, Q)$, there is
a zero $1/2 + i\gamma$ of the corresponding Epstein zeta-function
$\E$ with
\begin{equation}\label{sj-result}
|\gamma -T |\leq T^{5/11 + \varepsilon }.
\end{equation}
Thus the best answer to Question 1 for $\E$ is $H = T^{ 5/11 + \varepsilon }$, at the time
of writing this manuscript. However, in analogy with the $L$-functions associated to a cusp 
forms (mentioned in section \ref{intro}), the next big challenge seems to be to reduce the 
exponent of $T$ from $5/11$ to $1/3$ !
 
In this article, we prove the following 

\begin{theorem}\label{eps-gap}
Suppose $V$ satisfies
$$T^{\epsilon}\ll V \ll T^{1/2-\epsilon}$$
and $R:=R(V)$ denote the number of gaps of length at least $V$ between 
consecutive zeros of $\E(\half + it)$ in the interval $[0,T]$, then
\begin{equation}\label{eps}
R \ll T^{1+\epsilon} V^{-2}
\end{equation}
for sufficiently large $T$, where the constant in $\ll$ may depend only on $\Delta$ and  $\epsilon$.
\end{theorem}

\begin{rem}
The method of the proof closely follows the paper by Ivi\'{c} and Jutila \cite{ivic-jutila}.
Using their approach to find an analogue of 
\eqref{thmA2} and \eqref{thmB2} for $\E(s)$ does not seem to
yield anything better than \eqref{eps}. 
Note that an analogue of \eqref{thmA2} or \eqref{thmB2} would immediately
improve the exponent of $T$ to $1/3$  in \eqref{sj-result}.

\end{rem}

\section{Some preliminary lemmas}\label{lemmas}

Ramachandra showed (see \cite{Kram}, Chapter II ) that the first power mean 
of a generalized Dirichlet series
satisfying certain conditions can not be too small. The following Lemma ( Theorem 3 of \cite{B} ) is a 
particular case of this general theorem, which
 is quite useful in  obtaining lower bounds of this type, 
even in short-intervals.

\begin{lem}\label{lem3}
Let $  B(s) = \sum_{n=1}^{\infty} b_n n^{-s} $ be any Dirichlet series satisfying the following conditions:
\begin{itemize}
\item[(i)] not all $b_n$'s are zero;
\item[(ii)] the function can be continued analytically in $ \sigma \geq a, \  |t| \geq t_0 $, and in this
region $ B(s) = O( (|t|+10)^A ).$
\end{itemize}
Then for every $\epsilon > 0, $ we have 
$$  \int_T^{T+H} | B(\sigma + it)| dt \gg H $$
for all $ H \geq (\log T)^{\epsilon}, \  T\geq T_0(\epsilon),$ and $\sigma > a $.
\end{lem}

\medskip
A suitable Dirichlet polynomial approximation is required to 
replace the divergent series $\E(1/2+it)$ in the course of the proof. We state the
following Lemma which is a direct adaptation of  
 Lemma 3 in \cite{jutila1}. 

\begin{lem}\label{appx}
Let $t\geq2$ and $t^2\ll X \ll t ^A $, where $A$ is an arbitrarily
large positive constant. Then we have
\begin{equation}
\begin{gathered}
\zeta_{Q}\left(1/2 + it\right)=\sum_{n \leq X} r_Q (n)
n^{-1/2 -it}
\\ +\,(\log2)^{-1}\sum_{X<n\leq2X} r_Q(n) \log (2X/n)
n^{-1/2 - it}\\ +\,(\log 2)^{-1} 2\pi \di^{-1/2}
\left(1/2 -it\right)^{-2}((2X)^{1/2 -it} -X ^{1/2 -it})+
O(tX^{-1/2}).
\end{gathered}
\end{equation}
\end{lem}

\medskip
The following result is Lemma 4.3 of Titchmarsh \cite{T}.

\begin{lem}\label{tit}
Let $F(x)$ and $G(x)$ be real functions, $G(x)/F(x) $ monotonic, and 
$F''(x)/G(x) \geq m > 0,$ or $ \leq -m < 0.$ Then
\begin{equation}\label{expint}
\left| \int_a^b G(x) e^{iF(x)} dx \right| \leq \frac{4}{m}
\end{equation}
\end{lem}

\medskip

\section{Summation and transformation formulae}\label{summation}

Let 
$$S=\sum_{n} \eta(n) r_Q (n) n^{-1/2 -it}$$
and
$e_k(x)=\exp(2\pi x/ k).$  
Introducing an
additive character $\pmod k$, we may write it formally as
\begin{equation}\label{summ1e}
S = \sum_{n} \eta(n) r_Q (n) n^{-1/2 -it} \e_k(nh)\cdot\e_k(-nh);
\end{equation}
as in \cite{Jut2}, the purpose of the extra exponential factor is
to damp the oscillations of the original exponential sum before an
application of the Voronoi summation. Thus $S$ is of the general
form
$$
 \sum_n r_Q (n) f(n) \e_k (hn),
$$
and a Voronoi summation formula for such sums was given in
\cite{Jut3}, Eq. (28). To state it, we need  some notation. The
Gauss sums related to the form $Q$ and additive characters are
$$
G_Q(k,h)=\sum_{x,y \pmod k} \e_k(hQ(x,y)),
$$
and it holds (see \cite{Sm}, Lemma 1)
\begin{equation}\label{summ2e}
|G_Q(k,h)|\leq(\di, k)k.
\end{equation}
Further, the summation formula involves an integral positive
definite quadratic form $Q^* (x,y)$ depending on $Q$ and $k$, and
the discriminant of $Q^*$ is at most $\di$ in absolute value.
Also, there
occurs an arithmetic function (corresponding to $\tilde{r}_{Q^*}(n)$
in \cite{Jut3}) of the form
\begin{equation}\label{summ3e}
\rho (n) = \rho (n; Q,h/k) = \sum_{Q^*(x,y)=n}\alpha(x,y),
\end{equation}
where $|\alpha(x,y)|\leq (\di, k) \leq \di$ and  $\alpha(x,y)$ depends only on
the classes of $x$ and $y$ $\pmod{\di}$ for given $Q$ and $h/k$.
 In this notation, a slightly simplified version of
the summation formula (28) in \cite{Jut3} can be stated as follows.

Let  $\chi(s)$ be as in the functional equation
$\zeta(s)=\chi(s)\zeta(1-s)$, thus $\chi(s)=2^s \pi^{s-1}
\sin(\frac{1}{2}\pi s)\Gamma(1-s)$.
The following Lemma is due to Jutila and Srinivas (Lemma 3.2 \cite{srini-jut})
\begin{lem}\label{transform}
Let $t$ be a large positive number, $r=h/k$ a positive rational
number with $(h,k)=1$, and suppose that the positive numbers $M_1$
and  $M_2$ satisfy
$$
\begin{gathered}
M_1 <\frac{t}{2\pi r}< M_2 \\
M_j =\frac{t}{2\pi r}
+(-1)^j m_j, \quad  m_1\asymp m_2,\\
1\leq k \ll {M_1}^{1/2-\delta_1},\\
t^{\delta_2} \mathrm{max}(t^{1/2} r^{-1}, hk)\ll m_1 \ll
{M_1}^{1-\delta_3}
\end{gathered}
$$
for some small positive constants $\delta_j$. Further, let
$$
U\gg r^{-1} t^{1/2+ \delta_4}
$$
and let $J$ be a fixed positive integer exceeding a certain bound
which depends on $\delta_4$.
Write
$$
M_j '= M_j + (-1)^{j-1} JU=\frac{t}{2\pi r} + (-1)^j m_j ' ,
$$
supposing  that $ m_j \asymp m_j '$, and let
\begin{equation}\label{njmj}
n_j =\di_0 h^2 {m_j}^2 {M_j}^{-1}, \quad n_j '=\di_0 h^2 (m_j ')^2
(M_j ')^{-1}.
\end{equation}
Then for a certain weight function $\eta \in C^{J-1}(\R)$ with
support $[\, M_1, M_2]$ and satisfying $\eta(x)=1$ for $x\in
[\,{M_1}',\, {M_2}'\,]$, we have
\begin{equation}\label{trans}
\begin{gathered}
\sum_{n=1} ^{\infty} \eta(n) r_Q (n) n^{-1/2 -it} = \biggl{\{}2\pi
\di^{-1/2} k^{-2} G_Q (k, -h) r^{-1/2} \\ +\, \pi^{1/4}i
(2hkt\di_0 )^{-1/4}(\di_0 /\di)^{1/2} \sum_{j=1}^2 (-1)^j
\sum_{n<n_j} w_j (n)\rho(n)  \\
\times \, \exp\left(2\pi i n\left(\frac{\bar{h} \bar{\di}
_0}{k}-\frac{1}{2hk\di _0}\right)
\right)n^{-1/4}\left(1+\frac{\pi n}{2hk\di _0 t}\right)^{-1/4} \\
\times \, \exp\left(i(-1)^{j-1}
\left(2t\phi\pfrac{\pi n}{2hkt\di_0} +\frac{\pi}{4}\right)\right)
\biggl{\}} r^{it}\chi\left(\frac{1}{2} +it\right)\\ +\, O(h^2
k^{-1} {m_1}^{1/2} t^{-3/2} U \log t),
\end{gathered}
\end{equation}
where
\begin{equation}\label{phi}
\phi(x)=\mathrm{arcsinh}(x^{1/2})+(x+x^2)^{1/2},
\end{equation}

$w_j(n)=1$ for $n<n_j'$ , $w_j  (n)\ll 1$ for $n\leq n_j$,
$w_j(y)$ and $w_j '(y)$ are piecewise continuous in the interval
$(n_j' , n_j)$ with at most $J-1$ jumps, and
$$
w_j '(y) \ll (n_j- n_j ') ^{-1} \quad \text{for } n_j ' <y < n_j
$$
whenever $w_j '(y)$ exists.
\end{lem}

\section{Proof of the theorem}\label{proof}

We
define the functions $f(s),\gamma(s)$ and $W(t)$ as
$$
f(s)=e^{\frac{1}{2} \pi i (\frac{1}{2}-s)}\left(
\frac{\sqrt{\di}}{2\pi}\right) ^s \Gamma(s) \E(s)=\gamma(s)\E(s)
$$
and
$$
W(t)=f\left(1/2 +it\right);
$$
the latter is an analogue of Hardy's function $Z(t)$ in the theory
of Riemann's zeta-function. The functional equation for $\E(s)$ 
implies that  $W(t)$ is real for real $t$. Thus the
zeros of $\E(s)$ on the critical line correspond to the real zeros of
$W(t)$. To prove the result we can obviously restrict ourselves to the 
zeros lying in $ [ T/2, T]$.
From now onwards $c_i$'s will always denote absolute constants.
Suppose now that $\tau$ and $\tau+U$ are the ordinates of two consecutive 
zeros of
$\E(1/2+it)$ and hence two consecutive zeros of $W(t)$ such that
\begin{equation}\label{tauT}
[\tau,\tau+U]\subset [T/2,T],\ \ \textrm{and} \ U\ge V, 
\end{equation} 
where $T^{\epsilon} \le V \le T^{5/11}$.
Further, let $L=\cc (\log T )^{\dd}$ and $G=VL^{-1}$. 
Then, for $t\in [\tau+V/4, \tau+ 3V/4]$,
$W(u)$ has no zero in the interval $[t-V/4,t+V/4]$ 
and so $W(u)$ does not change sign in this interval.
Let us define
$$I_1(t)=\int_{t-V/4}^{t+V/4}|W(u)|e^{-(t-u)^2/G^2} du$$
and
$$I_2(t)=\int_{t-V/4}^{t+V/4}W(u)e^{-(t-u)^2/G^2} du$$
Then for $t\in [\tau+V/4, \tau+ 3V/4]$,
$I_1(t)=|I_2(t)|$.
Let $\mathcal S$ be the set of $t\in [T/2+V/4, T-V/4]$ such that
$I_1(t)=|I_2(t)|$. 
Since the number of pairs of consecutive zeros of
$\zeta_Q(1/2 + i t)$, $\tau$ and $\tau+U$,  
such that \eqref{tauT} holds is $R$, 
the above discussion implies that the Lebesgue measure of
$\mathcal S$, $m(\mathcal S) \gg RV $.
Thus
\begin{equation}\label{RV}
RV \ll m(\mathcal S)=\int_{t\in \mathcal S} dt 
=\int_{t\in \mathcal S}\frac{|I_2(t)|^2}{I_1(t)^2} dt.
\end{equation}
Therefore, to obtain an upper bound on $R$,  
we have to find a  lower bound for $I_1(t)$ and an
upper bound for $I_2(t)$. 

\medskip
The lower bound 
\begin{equation}\label{lowerbd}
I_1(t)\gg G.
\end{equation}
follows from Lemma \ref{lem3}.

\medskip

From now onwards we shall dwell on the upper bound estimation of $I_2(t)$.
First we define some notations; let
$P(u)= {u\sqrt{\Delta}}/{2\pi}$, 
for $J$ sufficiently large, we choose a   
smoothing function $\eta \in C^{J-1}(\R)$ satisfying the following
\begin{equation}\label{eta}
\eta(x) = \left\{ 
\begin{array}{l l}
  1 & \quad \mbox{if} ~~ x\in [P(T)-Y,P(T)+Y]\\
  0 & \quad \mbox{if} ~~ x>  P(T)+ 2Y ~~\textrm{or}~~ x < P(T)-2Y \\
\end{array} \right.
\end{equation}
The parameter $Y$ is defined by
\begin{equation}\label{kv}
VY=T^{1+\epsilon}.
\end{equation}

We shall use $e^x$ and $\exp(x)$ to denote the same function. 
Now the integral  $I_2(t)$ is written as  
\begin{align}\label{I2t1}
I_2(t)=\int_{-V/4}^{V/4} \, &  \gamma({1}/{2} +i(t+u))\, 
   \E({1}/{2} + i(t+u)) \, e^{-(u/G)^2} du.
\end{align}

We substitute the expression for $\E(1/2 + i(t+u))$ given by Lemma \ref{appx} 
with $X=t^3$ in \eqref{I2t1} to obtain
\begin{align}\label{I2t2}
& I_2(t) =   \sum_{\substack{|n-P(t)|>2Y \\ n\le t^3 }}r_Q(n)n^{-1/2}
             \int_{-V/4}^{V/4}\gamma(1/2+i(t+u))n^{-i(t+u)}e^{-(u/G)^2}du\\
&  \quad +   \int_{-V/4}^{V/4}\gamma(1/2+i(t+u))
             \left(\sum_{|n-P(t)|\le 2Y}r_Q(n)n^{-1/2-i(t+u)}\right) 
             e^{-(u/G)^2}du \nonumber\\
&  \quad +   (\log 2)^{-1}\sum_{t^3<n\le 2t^3}r_Q(n)\log (2 t^3/n)n^{-1/2}\nonumber\\
&  \quad \times  \int_{-V/4}^{V/4}\gamma(1/2+i(t+u))n^{-i(t+u)}e^{-(u/G)^2}du
             + O(G t^{-1/2})\nonumber\\
&  \quad =   S_1+S_2+S_3+o(G).\nonumber
\end{align}

We expand the sum within brackets in $S_2$ and use \eqref{eta} to get
\begin{align}\label{I2t3}
S_2  &=  \int_{-V/4}^{V/4}\gamma(1/2+i(t+u)) \Bigg(\sum_{|n-P(t)|\le2Y}
         \eta(n)r_Q(n)n^{-1/2-i(t+u)} \\
&   +  \sum_{|n-P(t)|\le 2Y}(1-\eta(n))r_Q(n)n^{-1/2-i(t+u)}\Bigg) e^{-(u/G)^2}du\nonumber\\
&   =  \int_{-V/4}^{V/4}\gamma(1/2+i(t+u))\sum_{n=1}^{\infty}\eta(n)
    r_Q(n)n^{-1/2-i(t+u)}e^{-(u/G)^2}du\nonumber\\
&   +  \sum_{Y<|n-P(t)|\le 2Y}(1-\eta(n))r_Q(n)n^{-1/2} \nonumber\\
&   \times \int_{-V/4}^{V/4} \gamma(1/2+i(t+u)) n^{-i(t+u)}e^{-(u/G)^2}du\nonumber\\
&   =  S_4+S_5.\nonumber
\end{align}

\medskip

We now estimate $S_1,S_3$ and $S_5$. Note that in all these cases $|n-P(t)|>Y$ and $n \ll t^3$.
First we consider $n>P(t)+Y$. We view the integral
$$\int_{-V/4}^{V/4}\gamma(1/2+i(t+u))n^{-i(t+u)}e^{-(u/G)^2}du$$
appearing in  $S_1,S_3,S_5$,
as a complex integral over the rectangular contour with vertices
$\pm V/4$, $\pm V/4 -i G$. We then use the well known 
 Stirling's formula for the $\Gamma$-function which states
that in any fixed vertical strip 
$-\infty < \alpha \leq \sigma \leq \beta < \infty$,
$$\Gamma(\sigma+it)=(2\pi)^{1/2} t^{\sigma+it-1/2}e^{-\frac{\pi}{2} t + 
\frac{\pi}{2} i (\sigma-1/2)-i t} ( 1 + O(1/t) ) \quad \hbox{as} \ t \rightarrow \infty ,$$
in order to estimate $\gamma(1/2+i(t+u))$ along this contour.

Therefore, on the vertical sides, where $u=\pm V/4-iw$ with $w\in [0,G]$, we estimate the factors in the integrand as follows
\begin{align*}
& \gamma(1/2+i(t+u))  n^{-i(t+u)}  = \Delta^{1/4}\bigg(
 \frac{2\pi n}{\sqrt{\Delta}(t+\frac{V}{4})}\bigg)^{-w} \times\\
& \times \exp{\textstyle \left\lbrace 
i\left( \left(t+\frac{V}{4}\right)\log\left( \frac{ \sqrt{\Delta}}{ 2\pi n}
\left(t+\frac{V}{4}\right)\right) 
-\left(t+\frac{V}{4}\right)\right) \right\rbrace}\big(1+O(1/T)\big)\\
&\ll \bigg(\frac{n}{P(t)+\frac{\sqrt{\Delta}}{8\pi}V}\bigg)^{-w} 
\end{align*}
which is bounded. As for the other factor, we have
$$e^{-(u/G)^2}\ll e^{-(L/4)^2} = e^{- 4 \log T }.$$
Hence the value of the integral on the vertical sides becomes $O(G/T^4)$.

On the horizontal side in the  lower half plane, 
let $u=v-iG$ with $v\in [-\frac{V}{4},\frac{V}{4}]$. 
Here $e^{-(u/G)^2}$ is bounded and the other factor is estimated by
\begin{align*}
& \gamma(1/2+i(t+u)) n^{-i(t+u)}
\ll \exp \bigg( -G\log \bigg(\frac{n}{P(t)+\frac{\sqrt{\Delta}}{8\pi}V}\bigg) \bigg) \\
&\ll \exp\left(-\frac{G Y}{P(T)}\right)
= \exp(-c\frac{T^{\epsilon}}{L}) \ll \exp(-\log(T^4 L))
\end{align*}
where the constant $c=\sqrt{\Delta}/(2\pi)$. Hence the value of the integral on the horizontal side is also $O(G/T^4)$.

For $n<P(t)-Y$  we take the rectangular contour with vertices
$\pm V/4$, $\pm V/4 +i G$ in the upper half plane and argue in a similar way to estimate the value of the integral to be $O(G/T^4)$.

Finally we use the fact that $r_Q(n) \ll n^{\epsilon_1}$ for any $\epsilon_1 >0$ (see the discussion preceding \eqref{pf9}) and estimate $S_1,S_3$, and $S_5$ as follows
\begin{align}\label{Sj}
 S_j \ll G T^{-4} \sum_{n\leq 2t^3}r_Q(n)n^{-1/2} 
     \ll G T^{-4 + 3/2 +\epsilon_1} = o(G) \qquad j=1,3,5.
\end{align}

From equations \eqref{I2t2},\eqref{I2t3} and \eqref{Sj}, after a change of variable, we conclude that
\begin{equation}\label{I2tfinal}
I_2(t) \ll \int_{t-V/4} ^{t+V/4} 
         |\sum_{n=1}^{\infty} \eta(n)r_Q(n)n^{-1/2-iu}|
          e^{-(t-u)^2/G^2}du +o(G).
\end{equation}

\medskip
Now we need to use the transformation formula \eqref{trans} of Lemma \ref{transform} to 
transform the sum in \eqref{I2tfinal}.
For that we first choose a rational approximation $r$ to $1/\sqrt{\Delta}$ as follows.
Since $\sqrt{\Delta}$ is a quadratic irrational, we can choose $r=h/k$ with $(h,k)=1$, 
such that $h^{-2} \ll |\sqrt{\Delta}-k/h| \ll h^{-2}$ and $h \asymp \sqrt{T/Y}$. 
We also observe that $k \asymp h$ as $h/k = r \asymp 1/\sqrt{\Delta},$ which is a constant. Let 
$M_j = P(u)+ (-1)^j 2Y$, $j=1,2$. Then 
$$
m_j = (-1)^j \frac{u}{2\pi} \left(\sqrt{\Delta}-\frac{1}{r}\right)+2Y 
\asymp \max(T h^{-2},Y) \ll Y.
$$
Also we choose
$$M_j' = P(u)+ (-1)^j Y, \ \ m_j' = (-1)^j \frac{u}{2\pi} 
  \left(\sqrt{\Delta}-\frac{1}{r}\right)+Y \ll Y .$$
We observe that $M_j, M_j'\asymp T$.
It is easy to check that all conditions of the Lemma \ref{transform} are satisfied when 
$T^{1/2+\epsilon} \ll Y \ll T^{1-\epsilon}$ or in other words when
$T^{\epsilon} \ll V \ll T^{1/2-\epsilon}$. 
Hence by Lemma \ref{transform}, there exists a smoothing function $\eta$ satisfying
\eqref{eta} for which \eqref{trans} holds.
Finally we also get the following estimate for $n_j(u)$, 
\begin{equation}\label{pf_nj}
 n_j(u) \ll \frac{h^2 m_j^2}{M_j} \ll Y.
\end{equation}
Now we replace the sum in \eqref{I2tfinal} by using \eqref{trans}. Then the 
leading term in \eqref{trans} is $\ll k^{-2}|G_Q(k,-h)| \ll k^{-1}$ 
which after integration in \eqref{I2tfinal} 
becomes $o(G)$.
Similarly, the error term in \eqref{trans} is $\ll h^2 k^{-1} Y^{1/2} T^{-3/2} U \log T
\ll T^{-1/2+\epsilon}$ which, after integration in \eqref{I2tfinal}, also becomes $o(G)$.

Thus we obtain the following upper bound on $I_2(t)$,
\begin{equation}\label{upperbd1}
I_2(t) \ll \int_{t- V/4}^{t+ V/4} |\Sigma(u)| e^{-(t-u)^2/ G^2} du +o(G),
\end{equation}
where $\Sigma(u)$ is as follows
\begin{align}\label{sigmau}
& \Sigma(u) =   (h k u)^{-1/4} r^{iu}\chi\left({1}/{2} +iu\right)
     \sum_{j=1}^2 (-1)^j \sum_{n<n_j(u)} w_j (n)\rho(n)n^{-1/4} \\
& \times  \, \exp(2\pi i C_1 n) \left(1+{C_2 n}/{u}\right)^{-1/4} \exp\left(i(-1)^{j-1}
\left(2u \phi(C_2 n /u) +\pi/4 \right)\right),\nonumber
\end{align}
where $C_1 = \bar{h} \bar{\di}_0/k-1/(2hk\di _0)$, and $C_2 = \pi / (2hk\di _0)$.

\medskip
 We now claim that for $t \in \mathcal S$ we can ignore the error term in \eqref{upperbd1}.
 This is because for these $t$ we have (by \eqref{lowerbd}) $G \ll I_1(t)=|I_2(t)|$. Hence
 we actually have
\begin{equation}\label{upperbd2}
I_2(t) \ll \int_{t- V/4}^{t+ V/4} |\Sigma(u)| e^{-(t-u)^2/ G^2} du.
\end{equation}

Using the lower bound \eqref{lowerbd} of $I_1(t)$ and the upper bound \eqref{upperbd2}
of $I_2(t)$ in \eqref{RV} we get
\begin{equation}\label{RV2}
RV \ll G^{-2}\int_{t\in \mathcal S}
\left|\int_{t- V/4}^{t+ V/4} |\Sigma(u)| e^{-(t-u)^2 /G^2} du\right|^2 dt.
\end{equation}

Applying the Cauchy-Schwarz inequality to the inner integral, \eqref{RV2} becomes
\begin{equation}\label{RV3}
RV \ll G^{-2} \int_{t\in \mathcal S}
    \left(\int \big(|\Sigma(u)|e^{-\half (t-u)^2 / G^2}\big)^2 du
 \int (e^{-\half (t-u)^2 / G^2})^2 du\right) dt,
\end{equation}
where the inner integrals are over the range $[t- V/4,t+ V/4]$. 
Note that
\begin{equation}\label{exp}
 \int_{t- V/4}^{t+ V/4} e^{-(t-u)^2 / G^2} du = 
 \int_{-V/4}^{V/4} e^{-x^2 / G^2} dx \leq \sqrt{\pi} G.
\end{equation}
Hence   \eqref{RV3}  and \eqref{exp}  yields
\begin{equation*}\label{pf4}
RV \ll G^{-1} \int_{t\in \mathcal S}\int_{t- V/4}^{t+ V/4}|\Sigma(u)|^2 e^{-(t-u)^2 / G^2}du.
\end{equation*}
By replacing $\mathcal S$ by the bigger set $[T/2+V/4,T-V/4]$ and interchanging the order of integration we get
\begin{equation}\label{pf5}
 R \, V \ll G^{-1} \int_{ T/2}^{T}  |\Sigma(u)|^2 
\int_{u- V/4}^{u+ V/4}
e^{-(t-u)^2 / G^2} dt \, du \ll \int_{T/2}^{T}  |\Sigma(u)|^2 du.
\end{equation}
Using \eqref{sigmau} in \eqref{pf5}, we have the following bound 
\begin{align}\label{pf6}
 R  & \ll V^{-1} (h k T)^{-1/2} \sum_{j=1}^2 \int_{ T/2}^{T}
 \bigg| \sum_{n<n_j(u)} w_j(n) \rho(n) n^{-{1}/{4}} \exp(2 \pi i C_1 n)   \\
  & \times \,  \left(1+C_2 n / u\right)^{-{1}/{4}}
 \exp\left(i(-1)^{j-1} 2u\phi(C_2 n / u)\right) 
\bigg|^2 du . \nonumber
\end{align}
Suppose that $n_j(T_0)=\max \{n_j(u):u\in [T/2,T]\}.$ By expanding the sum 
within the modulus in \eqref{pf6} we get
\begin{align}\label{R}
 R   &\ll  V^{-1} (h k)^{-1/2}T^{1/2} \sum_{j=1}^2 \sum_{n<n_j(T_0)} 
|w_j(n) \rho(n)|^{2} n^{-1/2}  \\
&  + \   V^{-1} (h k T)^{-1/2}  \sum_{j=1}^2 \sum_{\substack{m,n<n_j(T_0) \\ m\neq n}}
 |w_j(m)\rho(m)w_j(n)\rho(n)| (m n)^{-{1}/{4}} \nonumber \\
& \qquad\qquad\times \int ((1+C_2 {m}/{u})^{-{1}/{4}} (1+C_2 {n}/{u})^{-{1}/{4}} \nonumber\\
& \qquad\qquad\times \exp\left(i(-1)^{j-1} 2u\left(
 \phi\left(C_2 {m}/{u}\right)-\phi\left(C_2 {n}/{u}\right) \right)\right) du \nonumber \\
& =   \mathbf{\Sigma_1 + \Sigma_2}. \nonumber
\end{align}
In $\mathbf \Sigma_2$, the integral will be over an appropriate subinterval 
of $[T/2,T]$ depending on $m$ and $n$.

\medskip

For estimating $\mathbf{\Sigma_1}$ we will need an estimate on the mean-square of the coefficients $r_{Q^*}(n)$. Since $r_ {Q^*}(n)$ is a complex linear combination of the
Dirichlet coefficients of the Dedekind zeta-function $\zeta_K(s)$, we have,
for any $\epsilon>0$, 
$r_{Q^*}(n)\ll n^{\epsilon}$
for sufficiently large $n$ depending on $\epsilon$.
Therefore the following estimate trivially holds.
\begin{equation}\label{pf9}
 \sum_{n \leq x} r_{Q^*}^2(n) = O\left(x^{1+\epsilon}\right),
\end{equation}
for any $\epsilon > 0$. 
In fact, for certain quadratic forms much better results are known
(see \cite{murty-osburn}, \cite{osburn}).
Since $\rho(n) \ll r_{Q^*}(n)$ by \eqref{summ3e}, \eqref{pf9} implies that the
function $m(x)$ defined by
\begin{equation}\label{pf10}
 m(x) = \sum_{n \leq x} |w_j(n)\rho(n)|^2 \ll x^{1+\epsilon}
\end{equation}
for $x\le n_j(T_0)$.
Hence, by partial summation we get that the inner 
sum over $n$ in $\mathbf \Sigma_1$ is
$$
\begin{gathered}\label{pf11}
 \sum_{n<n_j} 
|w_j(n) \rho(n)|^2 n^{-1/2} = \int_0^{n_j} x^{-1/2} dm(x) = \\
= m(n_j)n_j^{-1/2} +
\half \int_0^{n_j} m(x)x^{-{3}/{2}} dx  \ll n_j^{1/2+\epsilon},
\end{gathered}
$$
where $n_j = n_j(T_0)$. Hence, 
\begin{equation}\label{sigma1}
\mathbf{ \Sigma_1} \ll V^{-1} (h k)^{-1/2}T^{1/2} Y^{1/2+\epsilon} 
\ll V^{-1} Y^{1+\epsilon} \ll T^{1+\epsilon}V^{-2}.
\end{equation}

\medskip

For estimating $\mathbf \Sigma_2$, we use Lemma \ref{tit} with 
$G(x)=((1+C_2 m /x)(1+C_2 n/x))^{-1/4}$ and 
$F(x) = 2 x(\phi(C_2 m/x)-\phi(C_2 n/x))$. From the definition \eqref{phi} of $\phi$ we get 
\begin{align*}
 F'(x) &= 2 \,\textrm{arcsinh}\left(\left(C_2{m}/{x}\right)^{1/2}\right)
-2\,\textrm{arcsinh}\left(\left(C_2{n}/{x}\right)^{1/2}\right) \\
& + 4\left(\left(C_2 {m}/{x}\right)\left(1+C_2 {m}/{x}\right)\right)^{1/2}
-4\left(\left(C_2 {n}/{x}\right)\left(1+C_2 {n}/{x}\right)\right)^{1/2}.
\end{align*}

Therefore $G(x)/F'(x)$ is monotonic in the interval $[ T/2,T]$, after
 replacing $F(x)$ by $-F(x)$ if necessary. We also have 
$|F'(x)/G(x)| \gg (C_2/T)^{1/2} |\sqrt{m}-\sqrt{n}|$. Hence, by 
Lemma \ref{tit} we get that the integral in $\mathbf \Sigma_2$ is 
$\ll (T/C_2)^{1/2} |\sqrt{m}-\sqrt{n}|^{-1} 
\ll (h k T)^{1/2} |\sqrt{m}-\sqrt{n}|^{-1}$. 
Thus we estimate $\mathbf \Sigma_2$ as follows.
$$
\begin{gathered}\label{pf12}
\mathbf{ \Sigma_2} \ll V^{-1} \sum_{j} \sum_{\substack{m , n < n_j \\ m\neq n}}
 |w_j(m)\rho(m)w_j(n)\rho(n)| (m n)^{-{1}/{4}} |\sqrt{m}-\sqrt{n}|^{-1}  \\
\ll V^{-1} \sum_{j} \sum_{n < m < n_j}
 |w_j(n)\rho(n)w_j(m)\rho(m)| n ^{-{1}/{4}} m^{{1}/{4}}(m-n)^{-1} \\
\ll V^{-1} \sum_{j} \sum_{h < n_j} h^{-1} \sum_{n  < n_j - h}
 |w_j(n)\rho(n)w_j(n+h)\rho(n+h)| n^{-{1}/{4}} (n+h)^{{1}/{4}},
\end{gathered}
$$
 where $n_j=n_j(T_0)$. Using the Cauchy-Schwarz inequality, the inner sum 
over $n$ is bounded by
 $$
\begin{gathered}
 \ll \bigg(\sum_{n < n_j} |w_j(n)\rho(n)|^2  n^{-1/2}\bigg)^{1/2}
 \bigg(\sum_{n+h < n_j}|w_j(n+h)\rho(n+h)|^2 (n+h)^{1/2}\bigg)^{1/2}.
\end{gathered}
$$
By using \eqref{pf10} and partial summation as was done in estimating
 $\mathbf \Sigma_1$, we conclude that the above is
$$
\begin{gathered}
 \ll \left(n_j^{1/2+\epsilon}n_j^{{3}/{2}+\epsilon}\right)^{\half} = n_j^{1+\epsilon}.
\end{gathered}
$$
Hence, we combine this with $\sum_{h < x} h^{-1} \ll \log x$ to 
finally get the following estimate
\begin{equation}\label{sigma2}
\mathbf{ \Sigma_2} \ll V^{-1} \sum_{j} n_j^{1+\epsilon} \log n_j 
          \ll V^{-1} Y^{1+\epsilon}
          \ll T^{1+\epsilon}V^{-2}. 
\end{equation}
Thus from \eqref{sigma1}, \eqref{sigma2} and \eqref{R}, we get
$$R \ll T^{1+\epsilon}V^{-2}.$$
This completes the proof of the theorem.

\end{document}